\documentclass[reqno,psamsfonts,11pt]{amsart}
\usepackage{cite}
\usepackage{multicol}
\usepackage{amscd}
\usepackage{graphicx}

\usepackage{amssymb}

\usepackage{latexsym}

\newcommand{\R}{\mathbb{R}}

\newtheorem{thm}{Theorem}[section]
\newtheorem{lem}[thm]{Lemma}

\newtheorem{pro}[thm]{Proposition}
\newtheorem{rmk}[thm]{Remark}

\begin{document}
\author[M. Abounouh, A. Atlas, O. Goubet]{}

\title[Boussinesq system]{Large time behavior of solutions to a dissipative Boussinesq system}

\email{abounouh@fstg-marrakech.ma,abdelghafour.atlas@u-picardie.fr \\  and olivier.goubet@u-picardie.fr}



\date{9 april 2006}

\maketitle

\centerline{\scshape  M. Abounouh}
 \medskip
 {\footnotesize \centerline{Universit\'e Cadi Ayyad, Facult\'e des
Sciences et Techniques }
 \centerline{Avenue Abdelkrim Khattabi, BP 618,
Marrakech, Maroc.   } }
  \medskip
\centerline{\scshape  A. Atlas, O. Goubet}
 \medskip
 {\footnotesize \centerline{ LAMFA CNRS
UMR 6140 }
  \centerline{Universit\'e de Picardie Jules Verne, } \centerline{33, rue Saint-Leu, 80039
Amiens, France.  } }
 \medskip

\begin{abstract}
In this article we consider the Boussinesq system supplemented with
some dissipation terms. These equations model the propagation of a
waterwave in shallow water. We prove the existence of a global
smooth attractor for the corresponding dynamical system.
\end{abstract}

\begin{center}
\section{Introduction}
\end{center}
This article is concerned with the long time behavior of the
solutions to a damped-forced Boussinesq system that read
    \begin{eqnarray}\label{FirstEquation}
      \left\{
            \begin{array}{lll}
             \eta_{t}+u_{x}+(\eta u)_x- \eta_{xx} &=&0,   \\
             u_{t}-u_{txx}- u_{xx}+\eta_{x}+u_x u &=&f.   \\
            \end{array}
      \right.
    \end{eqnarray}\\
Here we have an incompressible fluid on a channel. $u(t,x)$ is the
horizontal velocity at the top of the fluid, $\eta$ is the
fluctuation of the height of the fluid with respect to the rest
position that is $z=\eta(t,x)=0$, assuming that the bottom of the
channel is at $z=-1$. Observe that in our model we have to ensure
that $\eta(t,x)>-1$
$\forall t,x $.\\

Here $f(x)$ is an external force that does not depend on time and
the damping terms are respectively $-u_{xx},-\eta_{xx}$. In the
conservative case, that read
\begin{eqnarray}\label{boussinesq}
      \left\{
            \begin{array}{lll}
             \eta_{t}+u_{x}+(\eta u)_x &=&0,   \\
             u_{t}-u_{txx}+\eta_{x}+u_x u &=&0;   \\
            \end{array}
      \right.
    \end{eqnarray}\\

\noindent this system has been introduced by Boussinesq in 1877 to
model the fluctuation of a waterwave in shallow water. Other
well-known asymptotic models are Korteweg-de Vries equations and
Benjamin-Bona-Mahony equation, also known as the regularized long
wave equation. For these asymptotical models we would like to refer
to \cite{ChenBonaSaut, Whitham} and to the references therein.

In this article we are interested in the dissipative case. In the
case where $f=0$, the solutions converge to the equilibrium and the
issue is to find out the rate of convergence. Following the
pioneering work of Amick, Bona and Schonbek \cite{ABS}, this issue
has been addressed in the case where $x\in \mathbb{R}^D, D\geq 2$
using the famous Schonbek splitting method \cite{Sandco}. Here we
plan to study the dynamical system provided by (\ref{boussinesq})
into the framework of infinite dimensional dynamical system
\cite{Temam, Raugel, Hale}. Our main result states as follows

\begin{thm}
The dynamical system provided by (\ref{FirstEquation}) features a
compact global attractor into a suitable energy space. Moreover this
compact global attractor has finite fractal and Hausdorff dimension.
\end{thm}

\noindent This result compares with previous results obtained for
dissipative KdV equations \cite{Gh1, Gh2, G, GR} or dissipative BBM
equations \cite{Wang, Ak}. This article is organized as follows; in
the next section we introduce the mathematical framework that we
have chosen to study this dynamical system. In a third section we
address the initial value problem for the evolution equation. In a
fourth section we prove the existence of a smooth finite dimensional
attractor.

\section{Mathematical framework}
\subsection{Initial data}
For the sake of convenience, we are interested in considering
periodic boundary conditions. We now consider functions for
$x\in[0,1]$ that are $1$-periodic. We also
assume that $\int_0^1 f(x) dx=0 $ and $f\in L^2(0,1)$.\\
Introducing $w(t,x)=1+\eta(t,x)$, we rather use the following system

\begin{eqnarray}\label{newb}
      \left\{
            \begin{array}{lll}
             u_{t}-u_{txx}- u_{xx}+w_{x}+u_x u &=&f  \mbox{ \ \   for   }(x,t)\in(0,1)\times \R_{+}\, , \\
             w_{t}+(w u)_x- w_{xx} &=&0 \mbox{  \ \ for   }(x,t)\in(0,1)\times \R_{+}\, .\\
            \end{array}
      \right.
    \end{eqnarray}\\

The natural space for the velocity $u(t,x)$ is then
\begin{equation}
\dot{H}^1_{per}=\{v\in H^1_{per} / \int_0^1 v(x)dx=0\}.
\end{equation}\label{3}
\noindent We then assume that $u_0 \in \dot{H}^1_{per} $.\\
We now proceed to the assumptions on $w_0(x)$.\\
The first physical assumption is that
\begin{equation}\label{4}
\inf w_0(x)>0.
\end{equation}
This assumption ensures that the top of the fluid does not hit the
bottom of the channel. The second assumption is that
\begin{equation}\label{5}
\int_0^1 w_0=1,
\end{equation}
that describes that $w_0$ fluctuates around 1 the height at rest of
the fluid, and that the fluids has constant volume.\\
The third assumption is related to the very definition of the entropy
for convection equation, see \cite{Sandco}.\\
Introducing $Q(y)=y\ln y -y+1$ that is convex and non negative, we
assume that
\begin{equation}\label{6}
\int_0^1 Q(w_0(x))dx<+\infty.
\end{equation}

\begin{rmk}
Assume that $f=0$ here and that we are given regular enough solution
$(u,w)$ to $(2)$ such that $w>0$ $\forall x,t$ . Multiplying $(2)$
by $(u,1+\ln w)$  and summing the two resulting equations, we thus
obtain
\begin{equation}\label{7}
  \frac{d}{dt} \Big[\frac{1}{2}|u|^2_{H^1}+\int_0^1w \ln w
     \Big]+|u_x|^2_{L^2}+ \int_0^1\frac{w_x^2}{w}=0 ;
\end{equation}
then, using $\int_0^1w=1$
\begin{eqnarray*}
  \lim_{t\rightarrow \infty} \Big[\frac{1}{2}\int_0^1(u^2+u_x^2)+\int_0^1  Q(w)\Big]= 0
\end{eqnarray*}
and the fluid converges to the equilibrium $(u,w)=(0,1)$.
\end{rmk}

\subsection{Functional Analysis} 
 We set
 \begin{equation}
 \mathring{K} =\{w(x)>0;\quad \int_0^1 w(x)dx=1 \,\,\,and \,\,\,\int_0^1
Q(w(x))dx<\infty\};
 \end{equation}
 one may wonder which kind of topology we shall use in $ \mathring{K}$.
 First of all, we observe that $ \mathring{K}$ is a convex set (this is
 obvious since the map $y\rightarrow Q(y) $ is a convex function).
 Furthermore, $ \mathring{K}$ is related to the following Orlicz space
 (see \cite{Adams}). Introducing
 \begin{equation}
 H(y)=(1+y)\ln(1+y)-y,
\end{equation}
that is a convex function, we observe that $H$ and $Q$ are related
by the following inequalities.

\begin{lem} \label{10}
$\exists C_0>0$ such that $\forall y \geq 0$
\begin{equation}
  Q(y)-1\leq H(y)\leq C_0(Q(y)+1).
 \end{equation}
\end{lem}
Proof: On the one hand, if $y\in [0,1]$, then $Q(y)\leq 1-y \leq
1+H(y)$. If $y \geq 1$ then
$$Q(y)=y \ln y +1-y\leq (y+1) \ln(y+1)-y+1.$$
To establish the reverse inequality, we observe that $y\ln y
\sim_\infty (1+y) \ln (y+1)$, then, for large $y$'s , say $y \geq
R$,
$H(y)\leq 2(y\ln y -y+1)$.\\
For $y\in [0,R]$, $H(y)$ is bounded by $H(R)$.Then the proof of
the lemma is completed. $\square$\\

Therefore $w>0$ belongs to $ \mathring{K}$ iff $\int_0^1w =1$ and $w\in
L_H$, the Orlicz space, whose norm is defined by
\begin{equation}\label{11}
||w||_{L_H}=inf\{\lambda >0,\,\, \int_0^1 H(\frac{w}{\lambda})(x) dx
\leq 1\}
\end{equation}
We now endow $ \mathring{K} $ with the topology of $L_H$, that is given by
the distance
\begin{equation}\label{12}
  d(w_1,w_2)=||w_1-w_2||_{L_H}.
 \end{equation}

 \begin{rmk}
$ \mathring{K}$ is not a closed subset of $L_H$, but $K=\{w \geq
0; \,\, \int_0^1w=1 \,\,and \,\, \int_0^1Q(w)<\infty\} $ is.
Consider $w_n$ in $K$ such that $
||w_n-w||_{L_H}\rightarrow 0$.
 There exists $\lambda_n \rightarrow 0$ such that
\begin{equation}\label{13}
  \int_0^1 H(\frac{|w_n-w|}{\lambda_n})(x)dx\leq 2.
 \end{equation}
We shall use
\begin{equation}\label{lemma2}
H(y) \geq \frac{1}{2} (\sqrt{y+1}-1)^2.
\end{equation}
\noindent Then, setting $ v_n(x)=\lambda_n^{-1}|w_n(x)-w(x)|$,
  \begin{eqnarray*}
& & \int_0^1  \frac{|w_n-w|}{\lambda_n}=\int_0^1 v_n \\
&=&\int_0^1 (\sqrt{v_n+1}-1)^2+2\int_0^1 (\sqrt{v_n+1}-1)\\
 &\leq&4 \int_0^1 H(v_n)+4(\int_0^1
 (\sqrt{v_n+1}-1)^2)^{\frac{1}{2}}\leq 4(1+\sqrt{2}).
\end{eqnarray*}
Then $w_n \rightarrow w$ in $L^1$, and, up to a subsequence
extraction $w_n\rightarrow w $ a.e. $\square$
\end{rmk}

\section{The Initial Value Problem}
\subsection{Main Theorem}
\begin{thm}\label{1}
Consider the initial data $(u_0,w_0)$ in $H^1_{per}\times  \mathring{K}$.
Then there exists a unique solution for (\ref{newb})
$$(u(t),w(t))\in C(\R_+;H^1_{per})\times
C(\R_+; \mathring{K})$$ that satisfies moreover $\sqrt{w}$ $\in
L^2_{loc}(\R_+;H^1_{per})$.
\end{thm}

\begin{rmk}
In the theorem above, we would like to point out two important
facts:\\
$\ast$ If $\inf w_0(x)>0$, then, for any $t,x>0$, $w(t,x)>0$, that
is
physically relevant.\\
$\ast$ The dissipative Boussinesq system provides a smoothing effect
in the $w$ variable. In fact $\sqrt w $ belongs $
C(\R_+-\{0\};H^1_x)$
\end{rmk}

Proof of Theorem \ref{1}:\\

Existence: we first regularize the initial data
$(u_0^\epsilon,w_0^\epsilon)$ to construct smooth solutions
$(u^\epsilon(t),w^\epsilon(t))$in $C([0,T],H_{per}^1\times
H_{per}^1)$ for instance. We then prove some a priori estimates and
finally pass to limit. Since these methods are classical, we just
indicate below how to derive the a priori estimates, referring the
reader to \cite{At}, for details. For the sake of simplicity,
we drop the subscript $\varepsilon$
throughout the proof of the theorem.\\

First step: we first prove $w(t,x)$ is positive .\\
\noindent Consider $\alpha= \inf w_0>0$. Introduce
\begin{equation}\label{14b}
J(t,x)=\max(0, \alpha-w(t,x)).
\end{equation}
\noindent Using the Kato's inequality ( see \cite{Kato, Kavian} )
that reads (in the distribution's sense)
\begin{equation}\label{15}
w_{xx}sgn(w)\leq (|w|)_{xx} ,
\end{equation}
\noindent we thus obtain
\begin{equation}\label{16b}
J_t+(uJ)_x-J_{xx}\leq 0.
\end{equation}
Then, integrating in $x$, we have that
\begin{equation}
\int_0^1 J(t,x)dx \leq \int_0^1 J(0,x)dx=0 ,
\end{equation}
and $w(t,x)\geq \alpha$ a.e.\\
This result is related to parabolic Harnack's inequalities, see \cite{Evans}.\\

Second step: a priori estimate in $H_{per}^1\times  \mathring{K}$.\\
We begin with a technical lemma
\begin{lem}\label{lemma3} Consider $w> 0$ a smooth periodic function such
that $\int_0^1 w(x)dx =1$ then
\begin{equation}\label{18}
\int_0^1Q(w)dx \leq \Big(\int_0^1 \frac{w_x ^2}{w} \Big)^{\frac
{1}{2}}.
\end{equation}
\end{lem}

Proof: since $-\ln$ is convex, $Q(w)\leq w(w-1)-w+1.$ Then, using
once more $\int_0^1 w=1$,
\begin{equation}\label{19}
\int_0^1 Q(w)\leq \int_0^1 w^2-1.
\end{equation}
On the other hand, for any $\varphi$ smooth periodic function, then
\begin{equation}\label{20}
\varphi ^2(x)\leq 2 \Big(\int_0^1
\varphi^2\Big)^{\frac{1}{2}}\Big(\int_0^1
\dot{\varphi}^2\Big)^{\frac{1}{2}}+\int_0^1 \varphi^2.
\end{equation}
We apply this to $\varphi=\sqrt{w}$. Then
\begin{equation}\label{21}
\int_0^1 w^2 \leq ||w||_{L^\infty} \Big(\int_0^1 w \Big) \leq
\Big(\int_0^1 \frac{w_x^2}{w} \Big)^{\frac {1}{2}}+1
\end{equation}
this concludes the proof of the lemma $\square$.

We now proceed to the a priori estimates. We multiply (\ref{newb})
by $(u,1+\ln w)$. We integrate in $x$ over $[0,1]$ the resulting
equations and then sum to obtain
\begin{eqnarray}\label{22}
\frac{d}{dt}\Big[ \int_0^1\Big( w\ln
 w-w+1\Big)&+&\frac{1}{2}||u||^2_{H^1}\Big
]+\int(uw)_x \ln w \nonumber\\
& &+\int_0^1\frac {w_x^2}{w}+ \int_0^1 w_xu
+||u_x||^2_{L^2}=\int_0^1 fu.
\end{eqnarray}
Integrating by part, we observe that
\begin{equation}\label{23}
\int_0^1(uw)_x \ln w+\int_0^1 w_x u=0.
\end{equation}
On the other hand, using Young's and Poincar\'e-Wirtinger
inequalities, we obtain
\begin{eqnarray}\label{24}
-\int_0^1fu +||u_x||^2_{L^2}&&\geq
||u_x||^2_{L^2}-(\pi-\frac{1}{2})||u||^2_{L^2}-\frac{1}{4\pi-2}||f||^2_{L^2} \nonumber\\
&& \geq \frac{1}{2}||u||^2_{H^1}-\frac{1}{4\pi-2}||f||^2_{L^2}.
\end{eqnarray}
We thus obtain,
\begin{equation}\label{25}
\frac{d}{dt}\Big[ \int_0^1Q(w)+\frac{1}{2}||u||^2_{H^1}\Big
]+\frac{1}{2}||u||^2_{H^1} + \int_0^1\frac{w_x ^2}{w} \leq
\frac{1}{4\pi-2}||f||^2_{L^2}.
\end{equation}
We now infer from (\ref{19}) and (\ref{21}) that
\begin{equation}\label{26}
\int_0^1Q(w) \leq \Big(\int_0^1\frac{w_x ^2}{w}\Big) +\frac{1}{4}.
\end{equation}
We combine this inequality together with (\ref{22}), we integrate
with respect to $t$ (thanks to the Gronwall lemma), and thus obtain
\begin{equation}\label{27}
\int_0^1Q(w(t))+\frac{1}{2}||u(t)||^2_{H^1}\leq
\frac{1}{4}+\frac{1}{4\pi-2}||f||^2_{L^2}+e^{-t}(\int_0^1Q(w_0)+\frac{1}{2}||u_0||^2_{H^1}).
\end{equation}
Since $H$ is a convex function, for $\lambda \geq 1$ (observe
$H(0)=0$), and due to Lemma \ref{10}
\begin{equation} \label{28}
 \int_0^1H(\frac{w}{\lambda})dx\leq \frac{1}{\lambda}\int_0^1H(w)dx
\leq \frac{C_0}{\lambda}\Big[1+
\sup_{t\geq0}\int_0^1Q(w(t))\Big]\leq 1,
\end{equation}
for $\lambda$ large enough. Then $w(t)$ remains bounded in the
Orlicz space $L_H$. We then have established an a priori estimate
for $(u,w)$ in $L^\infty(R_+;H_{per}^1)\times
L^\infty(R_+;L_H)$.

\begin{rmk}
Actually (27) implies that $\sqrt{w}$ is a.e. in $t$ in $H^1_x$.
Since this Sobolev space is an algebra and since we can solve the
evolution equation under consideration with initial data $w$ in
$H^1_x$, this implies that for all $t>0$ $w$ is in $H^1_x$
(smoothing effect). We precise this fact below.
\end{rmk}

 \noindent Third step: smoothing effect; $\sqrt{w}$
belongs to $H^1$
for $t>0$.\\
We set $v=\sqrt{w}$ that solves
\begin{equation}\label{29}
 v_t-v_{xx}-\frac{v_x^2}{v}+\frac{1}{2}u_xv+uv_x=0.
\end{equation}

Multiply (\ref{29}) by $-v_{xx}$ and integrate. We then get
\begin{eqnarray}\label{29,5}
\frac{1}{2}\frac{d}{dt}||v_x||_{L^2}^2&+&||v_{xx}||_{L^2}^2+\frac{1}{3}\int
\frac{v_x^4}{v^2}=\frac{1}{2}\int u_xvv_{xx}+\int
uv_xv_{xx}\nonumber\\
&&\leq
\frac{1}{2}||v_{xx}||_{L^2}^2+c\Big[||v||^2_\infty||u_{x}||_{L^2}^2+||u||^2_\infty||v_{x}||_{L^2}^2\Big].
\end{eqnarray}
Then, since $||v||_{L^2}=1$,
\begin{equation}\label{30}
\frac{d}{dt}||v_x||_{L^2}^2\leq c_1||u||^2_{H^1}||v||^2_{H^1} \leq
c_1||u||^2_{H^1}(1+||v_x||_{L^2}^2).
\end{equation}
Due to Gronwall lemma and since $u$ is bounded in $H^1$ we then get
the $H^1$ bound on $v=\sqrt{w}$.

\begin{rmk}
Actually for $T<+\infty$ $\sqrt{w}$ is in $L^2(0,T; H^1_x)$ and then
$w$ in $L^1(0,T; H^1_x)$.
\end{rmk}

\noindent Fourth step: uniqueness

Consider two trajectories $(u_2,w_2)$ and $(u_1,w_1)$ that start
from the same initial data. Due to the previous estimates both
$(u_2,w_2)$ and $(u_1,w_1)$ remain bounded in $L^1(0,T;
L_x^\infty)$. We set $u=u_2-u_1,w=w_2-w_1$ that are solutions to

\begin{eqnarray}\label{uniq}
      \left\{
            \begin{array}{lll}
             u_{t}-u_{txx}- u_{xx}+w_{x}+\frac{1}{2}(u_2^2-u_1^2)_x=0,  \\
             w_{t}- w_{xx}+(u_2w_2-u_1w_1)_x=0.\\
            \end{array}
      \right.
    \end{eqnarray}\\
\noindent Then multiply these equations by $(u,w)$ and integrate the
resulting equation to obtain

$$ \frac{d}{dt}( ||u||_{H^1}+||w||_{L^2})\leq C(1+
\max(||u_2||_{L^\infty},||w_2||_{L^\infty},||u_1||_{L^\infty},||w_1||_{L^\infty})
( ||u||_{H^1}+||w||_{L^2}).$$ \noindent The results follows
promptly. $\square$

\section{The global attractor}

\subsection{Existence of the global attractor}

To begin with we state and prove
\begin{pro}\label{un}
The semigroup $S(t)$ defined on $\dot{H}_{per}^1\times  \mathring{K} $
possesses an absorbing set that is bounded in $\dot{H}_{per}^1\times
H^1$
\end{pro}

\noindent Proof: The existence of a bounded absorbing set in
$\dot{H}_{per}^1\times  \mathring{K} $ comes from the estimate (\ref{27})
of the previous section. Let $t_0$ be the entrance time into this
absorbing ball. Going back to (\ref{30}) and applying the Uniform
Gronwall Lemma (see Lemma III.1.1 in \cite{Temam}), we thus obtain that
for $t>0$, for some numerical constant $c$,

\begin{equation}\label{ugl}
t||(\sqrt{w})_x(t+t_0)||^2_{L^2}\leq
c(1+t)(1+||f||^2_{L^2}))\exp(c+c||f||^2_{L^2}).
\end{equation}

\noindent Therefore $\sqrt{w}$ is bounded for large times into
$H^1$. Since $H^1$ is an algebra, then $w$ is also bounded for large
times in $H^1$. $\square$

\begin{thm}
The semigroup $S(t)$ possesses a global attractor $\mathcal{A}$ in
$\dot{H}_{per}^1\times L_H $, that is a compact subset of $H^2\times
H^2$.
\end{thm}
Proof: we introduce the splitting $(u,w)=(u^1,w)+(u^2,0)$,\\
where $u^1$ satisfies
\begin{eqnarray}\label{32}
      \left\{
            \begin{array}{lll}
             u^1_{t}-u^1_{txx}- u^1_{xx}+w_{x}+u_x u =f   \\
             u^1(0)=0,
            \end{array}
      \right.
    \end{eqnarray}\\
 and $u^2$ is solution to
\begin{eqnarray}\label{33}
      \left\{
            \begin{array}{lll}
             u^2_{t}-u^2_{txx}- u^2_{xx} =0   \\
             u^2(0)=u_0.
            \end{array}
      \right.
    \end{eqnarray}\\
We now define the families $\{S_1(t)\}_{t\geq 0}$ and
$\{S_2(t)\}_{t\geq 0}$ of maps in $H^1 \times L_H$, where
$S_1(t)(u_0,w_0)=(u^1,w)$ and $S_1(t)(u_0,w_0)=(u^2,0)$. \\
First step : we prove that $u^1$ is bounded in $H^2$. For this we
multiply (\ref{32}) $-u^1_{xx}$ and integrate between 0 and 1 to
obtain
\begin{eqnarray*}
\frac{1}{2}\frac{d}{dt}||u^1_x||_{H^1}^2+||u^1_{xx}||_{L^2}^2=-\int
f u^1_{xx}+\int w_x u^1_{xx}+\int uu_xu^1_{xx},
\end{eqnarray*}
due to Young and Cauchy-Schwarz inequalities, then
\begin{equation}\label{34}
\frac{d}{dt}||u^1_x||_{H^1}^2+||u^1_{xx}||_{L^2}^2\leq
c\Big(||f||_{L^2}^2+
||w_x||_{L^2}^2+||u_x||_{L^2}^2||u||_{L^\infty}^2\Big),
\end{equation}
 due to Proposition \ref{un}, (\ref{27}), Gronwall and Poincar\'e
inequalities, we obtain that $u^1$ remains in a bounded set of $H^2$
for large times ($t>t_0$ the entrance time into the absorbing ball).\\
On the other hand, it is an exercice to prove that
\begin{equation}
u^2(t)\rightarrow 0 \,\ strongly\,\, in \,\,H^1 \,\,when\,\,t\rightarrow
\infty.
\end{equation}
Then $S_1(t)(u_0,w_0)$ is bounded in $H^2\times H^1$ then compact in
$H^1\times L_H$
 and $S_2(t)(u_0,w_0) \rightarrow 0$ in $H^1\times L_H$ uniformly
 on bounded sets.\\
Then from Theorem I.1.1 in [T] we have the existence of a global
attractor $\mathcal{A}$ in $H^1\times L_H$, that is moreover a
bounded set in $H^2\times H^1$.\\
We now prove that for a trajectory $(u,w)$ in the global attractor,
then $w$ remains bounded in $H^2$. For that purpose, multiply the
second equation in (\ref{newb}) by $w_{4x}$ and integrate by parts
to obtain

\begin{equation}\begin{split}
\frac{d}{dt}||w_{xx}||_{L^2}^2+||w_{xxx}||_{L^2}^2=\int_0^1u_{xx}ww_{3x}+
2\int_0^1u_x w_xw_{3x}-\frac{1}{2}\int_0^1u_xw_{2x}^2 \\\leq
c||u||^2_{H^2}||w||^2_{H^1}+\frac{1}{2}||w_{xxx}||_{L^2}^2.
\end{split}\end{equation} \noindent Then the results follows
promptly. It remains to prove that the global attractor, that is
bounded in $H^2\times H^2$, is in fact a compact subset of this
space. This can be performed by the Energy Equation Method of
\cite{MRW} that is a suitable adaptation of the famous J. Ball
argument. This is standard and will not be reproduced here; we refer
the reader to \cite{At} for details. $\square$

\subsection{Dimension of the attractor}

In this section we are going to prove that the global attractor
$\mathcal{A}$ has a finite dimension in $\mathcal{E}=\dot{H}^1\times
\{w\in L^2; \int_0^1 w=1\}$. $\mathcal{E}$ is an affine space whose
associated vector space is $E=\dot{H}^1\times \dot{L}^2$. To begin
with, we need a result on the differentiability of the semi-group
$S(t)$ on the global attractor. Consider the non-autonomous
linearized system
\begin{eqnarray}\label{53}
            \left\{
        \begin{array}{lll}
             v_{t}-v_{txx}- v_{xx}+h_{x}+(u v)_x &=&0  \\
             h_{t}+(uh+vw)_x- h_{xx} &=&0 \\
            \end{array}
        \right.
    \end{eqnarray}\\
where $(u(t),w(t))=S(t)(u_0,w_0)$, $(u_0,w_0)\in \mathcal{E}$, is a
trajectory solution of (\ref{newb}) and $(v_0,h_0)\in E$. Actually
the linear mapping
 $DS(t)(u_0,w_0)(v_0,h_0)=(v(t),h(t))$ is the uniform differential
 of $S(t)$ as stated below
 \begin{thm}
The non-autonomous PDE (\ref{53}) provides a well posed initial
value problem in $E$. Moreover for $T>0$, $(v_0,h_0)\in E$,
$(u_0,w_0)\in\mathcal{A}$, $t\leq T$ there exists a constant
$C=C(T)$ such that
\begin{equation}
||S(t)(u_0+v_0,w_0+h_0)-S(t)(u_0,w_0)-DS(t)(u_0,w_0)(v_0,h_0)||_E\leq
C(T)||(v_0,h_0)||_E^\delta
\end{equation}
where $1<\delta<2$.
\end{thm}

\noindent Proof: to prove that the initial value problem is
well-posed is standard and then omitted. Consider the solutions
$(u_1(t),w_1(t))=S(t)(u_0,w_0)$,
$(u_2(t),w_2(t))=S(t)(u_0+v_0,w_0+h_0)$ and
$(v(t),h(t))=\Big(DS(t)(u_0,w_0)\Big)(v_0,h_0)$. Then
$(p,q)=(u_2,w_2)-(u_1,w_1)-(v,h)$ satisfies the system

\begin{eqnarray}\label{55}
      \left\{
            \begin{array}{lll}
             p_{t}-p_{txx}-p_{xx}+q_{x}+(\frac{1}{2}v^2+vp+u_1p+\frac{1}{2}p^2)_x &=&0  \\
             q_{t}- q_{xx}+(pq+p h +v q+h v+ q u_1+w_1p)_x &=&0 \\
            \end{array}
      \right.
    \end{eqnarray}

We shall use in the sequel that $\int_0^1 p=\int_0^1 q=0$ and then
$||p||_{H^1}$ and $||p_x||_{L^2}$ define equivalent norms. Multiply
(\ref{55}) by $(p,q)$ and integrate to obtain (due to
straightforward computations)

\begin{equation*}\begin{split}
\frac{1}{2}\frac{d}{dt}\Big[||q||_{L^2}^2+||p||_{H^1}^2\Big]+||q_x||_{L^2}^2
 +||p_x||_{L^2}^2 =\\ -\int q_x p+\frac{1}{2}\int (p+v)^2p_x +\int u_1pp_x
+   \int (pq+p h+ vq + v h+w_1p+ q u_1) q_x \\ \leq
(||v||_{H^1}+||h||_{L^2})||v||_{H^1}\Big[ ||q_x||_{L^2}
+||p_x||_{L^2}\Big]+\\
(1+||u_1||_{H^1}+||v||_{H^1}+||w_1||_{L^2}+||h||_{L^2})\Big[
||q_x||_{L^2}^2 +||p_x||_{L^2}^2\Big]+
\Big[||p||_{L^2}||q_x||_{L^2}^2 \Big].
\end{split}\end{equation*}

\noindent We thus obtain, using the bounds on the attractor and the
local in time bounds on $(v,h)$
\begin{equation}
\frac{d}{dt}\Big[||(p,q)||^2_E\Big]\leq K_1(T)||(p,q)||^2_E +
K_2(T)||(v_0,h_0)||^4_E+K_3(T)||(p,q)||^4_E.
\end{equation}

Consider a given interval of time $[0,T]$. Set
$\varepsilon^2=K_2(T)||(v_0,h_0)||^4_E$ that is small. Then
$\phi(t)=\exp(-tK_1(T))||(p,q)||^2_E(t)$ satisfies the ODE

\begin{equation}
\dot{\phi}\leq K \phi^2+\varepsilon^2,
\end{equation}
\noindent supplemented with $\phi(0)=0$. Then $E(t)\leq
2\varepsilon$ if $\varepsilon$ is small enough. $\square$\\

\noindent We now give the main result of this section
\begin{thm}
The fractal and Hausdorf dimension in $\mathcal{E}$ of the attractor
$\mathcal{A}$ are finite .
\end{thm}

\noindent Proof: set $\xi=(u,w)$, $\beta=(v,h)$. Now we study the
operators $DS(t)\xi_0$ that contracts the m-dimensional volumes in
$\mathcal{E}$. Let $\beta_0^1,...,\beta_0^m $ in $E$. We study the
following quantities
\begin{equation}
G_m=||\beta^1(t)\wedge ... \wedge \beta^m||_E^2=det_{1\leq i,j \leq
m}\Big(\beta^i(t),\beta^j(t)\Big)_E,
\end{equation}
where $\beta^i(t)=(DS(t)\xi_0)\beta^i_0$. The Gram determinant $G_m$
represents the volume of m-dimensional polyhedron defined by the
vectors $\beta^1(t),...,\beta^m(t) $. We will show that for
sufficiently large $m$ this determinant decays exponentially as
$t\rightarrow
\infty$. \\
We consider $\beta(t)=(DS(t)\xi_0)\beta_0$ solution of (\ref{53}),
we multiply by $\beta=(v,h)$ and integrate to obtain
\begin{equation}\label{58}
\frac{1}{2}\frac{d}{dt}||\beta||_E^2+||v_x||_{L^2}^2+||h_x||_{L^2}^2=\\
\int_0^1 uhh_x+\int_0^1 vwh_x+\int_0^1 uvv_x+\int_0^1 hv_x.
\end{equation}
\noindent Recall that $(u,w)$ is a trajectory that belongs to the
global attractor. Introduce $M=c(1+||f||^2_{L^2})$ that is the $H^1$
bound for $u$ in the attractor (see (\ref{27} )). We do not want to
use estimates that involve $H^2$ norms of $u$ as (\ref{ugl}).
\\
 We bound the
fourth term in the r.h.s of (\ref{58}) by
$\frac{1}{4}||h_x||_{L^2}^2+||v||_{L^2}^2$. The third term can be
bounded as follows
\begin{eqnarray*}\label{kk}
|\int_0^1 uvv_x|\leq
\frac{1}{2}||u||_{H^1}||v||_{L^2}||v||_{L^\infty}\leq \\
c||u||_{H^1}||v||^{3/2}_{L^2}||v_x||^{1/2}_{L^2}\leq
\frac{1}{4}||v_x||^2_{L^2}+cM^{4/3}||v||_{L^2}^2.
\end{eqnarray*}
We now proceed to the first term as follows
\begin{eqnarray*}\label{kkk}
|\int_0^1 uhh_x|\leq
\frac{1}{2}||u||_{H^1}||h||^2_{L^4}\leq \\
c||u||_{H^1}||h||_{H^{-1}}^{3/4}||h_x||^{5/4}_{L^2}\leq
\frac{1}{4}||h_x||^2_{L^2}+cM^{8/3}||h||_{H^{-1}}^{2}.
\end{eqnarray*}

\noindent For the second term, we have

\begin{eqnarray*}\label{kkkj}
|\int_0^1 vwh_x|\leq
\frac{1}{4}||h_x||^2_{L^2}+||w||^2_{L^2}||v||^2_{L^\infty}\leq \\
\frac{1}{4}||h_x||^2_{L^2}+\frac{1}{8}||v_x||^2_{L^2}+c||w||^4_{L^2}||v||^2_{L^2}.
\end{eqnarray*}
\noindent To go further, we need a new estimate on $w$ that reads
\begin{lem}
For any $(u,w)$ in $\mathcal{A}$, then $||w(t)||_{L^2}\leq
c(1+M^2)$.
\end{lem}
\noindent Proof: for a given trajectory in the attractor multiply
the second equation by $w$ and integrate to obtain
\begin{equation}\label{kkkjj}
\frac{d}{dt}||w||^2_{L^2}+2||w_x||^2_{L^2}=2\int_0^1 uww_x\leq
||w_x||^2_{L^2}+||u||^2_{L^\infty}||w||_{L^\infty};
\end{equation}
\noindent here we have used that $\int_0^1 w=1$. We then infer from
(\ref{kkkjj}) that, using Poincar\'e-Wirtinger inequality,
\begin{equation}
\int_0^1(w-1)^2=||w||^2_{L^2}-1\leq ||w_x||^2_{L^2},
\end{equation}
\noindent that
\begin{equation}
\frac{d}{dt}||w||^2_{L^2}+\frac{1}{4}||w||^2_{L^2}\leq c(1+M^4).
\end{equation}
\noindent Then the classical Gronwall lemma leads to the result.
$\square$

\noindent We then have
\begin{equation}\label{73}
\frac{1}{2}\frac{d}{dt}||\beta||_E^2+||\beta||_E^2=c\left((1+M^{8})
||\beta||_{L^2\times \dot{H}^{-1}}^2\right).
\end{equation}

We introduce the Gram determinant
$$G_m(t)=det_{1\leq i,j \leq m} \Big ( \Lambda( \beta^i(t) , \beta^j(t) ) \Big
)_E,$$ \noindent where
$\Lambda(a,b)=\frac{||a+b||_E^2-||a-b||_E^2}{4}$, and that
represents the $m$-dimensional volume. Then we can proceed as in
\cite{Temam, Gh1} to establish that
\begin{equation}
\frac{dG_m}{dt}+mG_m \leq c(1+M^{8})\left(\sum_{l=1}^m
\max_{A\subset \R^m,dimA=l}\min_{v\in A,v\neq
0}\frac{||\beta||_{L^2\times
\dot{H}^{-1}}^2}{||\beta||_E^2}\right)G_m.
\end{equation}
\noindent Since the eigenvalues of the Laplace periodic operator are
$4\pi^2k^2$ each of multiplicity $2$, then
\begin{equation}
\sum_{l=1}^m \max_{A\subset \R^m,dimA=l}\min_{v\in A,v\neq
0}\frac{||\beta||_{L^2\times \dot{H}^{-1}}^2}{||\beta||_E^2} \sim
2\pi^2 \sum_{k=1}^{m/2} (2\pi k)^{-2}\leq \frac{1}{12}.
\end{equation}
\noindent Therefore for $m\geq c(1+M^8)$ the $m$-dimensional volume
$G_m$ decays and the attractor has finite dimension. $\square$


 \vskip 1cm

 \centerline{AKNOWLEDGEMENT}
This work was performed under the research program "Dissipation et
Asymptotique" supported by CNRS (France) and CNRST (Marocco). The
third author would like to thank L. Dupaigne for pointing out to him
the reference \cite{Evans}. Thanks also to E. Zahrouni for helpful
remarks.

\vfill\break


\begin{thebibliography}{99}

\bibitem{Adams} R. A. Adams, {\sl Sobolev spaces}, Academic Press, New York, 1975.

\bibitem{Ak} N. Akroune, {\sl Comportement asymptotique de certaines \'equations faiblement amorties}, Phd (2000).

\bibitem{At} A. Atlas, {\sl PhD. Thesis}, Universit\'e de Picardie Jules Verne, 2006.

\bibitem{ABS}  C.~J. Amick, J.~L. Bona, and M.~E. Schonbek,
{\em Decay of solutions of
  some nonlinear wave equations}, J. Differential Equations, 81 (1989),
  pp.~1--49

\bibitem{ChenBonaSaut} M. Chen, J. Bona, JC. Saut  {\sl Boussinesq equations and other
systems for small-amplitude long waves in nonlinear dispersive media. I: Derivation and linear theory},   J. Nonlinear Sci.
{\bf 12}  (2002),  no. 4, 283--318.

\bibitem{Evans} L. C. Evans {\sl Partial differential equations}, Graduate Studies in Mathematics, 19. American Mathematical
Society, Providence, RI, 1998.

\bibitem{Gh1} J-M. Ghidaglia, {\sl Weakly damped forced Korteweg-de Vries
equations behave as a finite dimensional dynamical system in the
long time},  J. Diff. Eq., 74, pp 369-390, (1988).\par

\bibitem{Gh2} J-M. Ghidaglia,  {\sl A note on the strong convergence
towards attractors for damped forced KdV equations},  J. Diff. Eq. {\bf 110}, 356-359,
 (1994).

\bibitem{G} O. Goubet, {\sl Asymptotic smoothing effect for weakly damped forced
Korteweg-de Vries equations},  Discrete Contin. Dynam. Systems  {\bf 6}  (2000),
no. 3, 625--644.

\bibitem{GR} O. Goubet, R. Rosa {\sl Asymptotic smoothing and the global attractor
 of a weakly damped KdV equation on the real line },  J. Differential Equations  {\bf185}
   (2002),  no. 1, 25--53.


\bibitem{Hale} J. Hale, {\it Asymptotic behavior of Dissipative Systems},
Math. surveys and Monographs, vol 25, AMS, Providence, 1988.

\bibitem{Kato} T. Kato, {\sl  Schr\"odinger operators with singular potentials},  Israel J. Math.
{\bf 13}  (1972), 135--148.

\bibitem{Kavian} O. Kavian, {\sl Introduction \`a la th\'eorie des points critiques et applications
aux probl\`emes elliptiques}, Math\'ematiques et Applications vol.
{\bf13}, Springer, Paris, 1993.

\bibitem{MRW} I. Moise, R. Rosa and X. Wang, {\sl Attractors for non-compact semigroups
via energy equations},  Nonlinearity  {\bf 11},  (1998),  no. 5, 1369--1393.

\bibitem{Raugel}  G. Raugel, {\sl Global attractors in partial differential equations}.
Handbook of dynamical systems, Vol. 2,  885--982, North-Holland,
Amsterdam, 2002.

\bibitem{Sandco} M. E. Schonbek, {\sl Existence of solutions for the Boussinesq system of
 equations},  J. Diff. Eq. {\bf 42}, no. 3, 325--352,
 (1981).

\bibitem{Temam} R. Temam, {\it Infinite Dimensional Dynamical
Systems in Mechanics and Physics,\/} Springer-Verlag, Second
Edition, 1997. \par

\bibitem{Wang} B. Wang, strong attractors for the BBM equation, {\it
Appl. Math. Letters}, vol 10, n 2, pp 23-28, 1997.

\bibitem{Whitham} G. B. Whitham, {\sl Linear and nonlinear waves}, Wiley, New York, 1999,
reprint of the 1974 original, a Wiley-Interscience Publication.


\end{thebibliography}
\end{document}